\documentclass[12pt]{article}

\usepackage{amsmath}
\usepackage{amsthm}
\usepackage{amsfonts}
\usepackage{enumerate}
\usepackage[numbers,sort&compress]{natbib}
\usepackage{pstricks,pst-node}

\theoremstyle{plain}
\newtheorem{thm}{Theorem}
\newtheorem{lem}{Lemma}
\newtheorem{prop}[lem]{Proposition}
\newtheorem{cor}[thm]{Corollary}

\theoremstyle{definition}
\newtheorem{definition}{Definition}

\newcommand{\vertex}[2]{\cnode*(#2){4pt}{#1}}
\newcommand{\lbvertex}[2]{\cnode*(#2){6pt}{#1}\rput(#2){\white\tiny\textbf{\textsf{#1}}}}
\newcommand{\edge}[2]{\ncline[nodesep=0pt]{-}{#1}{#2}}

\title{Graphs with Given Degree Sequence and Maximal Spectral Radius}

\author{T\"urker B{\i}y{\i}ko\u{g}lu \\
  \small Department of Mathematics \\[-0.8ex]
  \small I\c{s}{\i}k University \\[-0.8ex]
  \small \c{S}ile 34980, Istanbul, Turkey \\[-0.8ex]
  \small \texttt{turker.biyikoglu@isikun.edu.tr} \\[2ex]
  Josef Leydold \\
  \small Department of Statistics and Mathematics \\[-0.8ex]
  \small University of Economics and Business Administration \\[-0.8ex]
  \small Augasse 2-6, A-1090 Wien, Austria \\[-0.8ex]
  \small \texttt{josef.leydold@wu-wien.ac.at}
}

\date{
  \small Mathematics Subject Classification: 
  05C35, 05C75, 05C05}

\begin{document}

\maketitle

\begin{abstract}
  We describe the structure of those graphs that have largest spectral
  radius in the class of all connected graphs with a given degree
  sequence. We show that in such a graph the degree sequence is
  non-increasing with respect to an ordering of the vertices induced
  by breadth-first search.
  For trees the resulting structure is uniquely determined up to
  isomorphism. We also show that the largest spectral radius in such
  classes of trees is strictly monotone with respect to majorization.
  \\[1ex]
  \emph{Keywords:}
  adjacency matrix, eigenvectors, spectral radius, degree sequence,
  Perron vector, tree, majorization
\end{abstract}


\section{Introduction}

Let $G(V,E)$ be a simple finite undirected graph with vertex set $V(G)$
and edge set $E(G)$. The eigenvalue of $G$ are the eigenvalues of the
adjacency matrix $A(G)$. The spectral radius of $G$ is the largest
eigenvalue of $A(G)$, also called the \emph{index} of the graph.
When $G$ is connected, $A(G)$ is
irreducible and by the Perron-Frobenius Theorem 
(see e.g.\ \citep{Horn:1990a})
the largest eigenvalue $\lambda(G)$ of $G$ is
simple and there is a unique positive unit eigenvector. We refer
to such an eigenvector $f$ as the \emph{Perron vector} of $G$.

There exists a vast literature that provides upper and lower bounds on
the largest eigenvalue of $G$ given some information about the graph,
for previous results see \citep{Cvetkovic;Rowlinson:1990a}.
Many recent results use the maximum, minimum or average degrees, e.g.,
\citep{Rojo:2006a, Liu;Shen;Wang:2007a}.
Some new results are based on the entire degree sequence, e.g.,
\citep{Shu;Wu:2004a}.

The goal of this article is slightly shifted. We want to characterize
connected graphs $G$ that have greatest spectral radius in the class
of all graphs with a given degree sequence. We show that in such
a graph the degree sequence is non-increasing with respect to an
ordering of the vertices induced by breadth-first search.
(Recently similar results have been shown for the special cases of 
caterpillars \citep{Simic;Marzi;Berlardo:2007a} and
cycles with spikes \citep{Belardo;Marzi;Simic:2006a}.)
We also show that the greatest maximum eigenvalue in such classes of trees 
is strictly monotone with respect to some partial ordering of degree
sequences.
The results are related to the (partly open) problem of finding
connected graphs of maximal spectral radius with given number of
vertices and edges (but arbitrary degree
sequences). \citet{Brualdi;Solheid:1986a} have shown 
that such graphs have stepwise adjacency matrix. We refer the reader
to \citep[Sect.~3.5]{Cvetkovic;Rowlinson;Simic:1997a} for details and
further discussion of this and related problems.

The paper is organized as follows:
The results of this paper are stated in
Section~\ref{sec:results}. In Section~\ref{sec:proofs} we prove these
theorems by means of a technique of rearranging graphs which has been
developed in \citep{Biyikoglu;Leydold:2006a} for the problem of
minimizing the first Dirichlet eigenvalue within a class of trees. 
Indeed, we will discuss the close relationship between this problem
and the problem of finding trees with greatest maximum eigenvalue in
Section~\ref{sec:remarks}.


\section{Degree Sequences and Largest Eigenvalue}
\label{sec:results}

Let $d(v)$ denote the degree of vertex $v$. We call a vertex $v$ with
$d(v)=1$ a \emph{pendant vertex} of the graph (and \emph{leaf} in case
of a tree). In the following $n$ denotes the total number of vertices,
i.e., $n=|V|$. A sequence $\pi=(d_0,\ldots,d_{n-1})$ of nonnegative
integers is called \emph{degree sequence} if there exists a graph $G$
with $n$ vertices for which $d_0,\ldots,d_{n-1}$ are the degrees of
its vertices, see \citet{Melnikov;etal:1994a} for relevant information.
In the entire article we enumerate the degrees in non-increasing
order. 

We introduce the following class for which we can provide optimal
results for the greatest maximum eigenvalue.
\begin{equation*}
  \mathcal{C}_\pi 
  = \{\text{$G$ is a connected graph with degree sequence $\pi$}\}\,.
\end{equation*}

For the characterization of graphs that have greatest maximum
eigenvalue among all graphs in $\mathcal{C}_\pi$ we introduce an
ordering of the vertices $v_0,\ldots,v_{n-1}$ of a graph by means of
breadth-first search:
Select a vertex $v_0\in G$ and create a sorted list of vertices
beginning with $v_0$; append all neighbors $v_1,\ldots,v_{d(v_0)}$ of
$v_0$ sorted by decreasing degrees; then append all neighbors of $v_1$
that are not already in this list; continue recursively with
$v_2,v_3,\ldots$ until all vertices of $G$ are processed.
In this way we build layers where each vertex $v$ in layer $i$ has
distance $i$ from root $v_0$ which we call its \emph{height}
$h(v)=\mathrm{dist}(v,v_0)$. Moreover, $v$ is adjacent to some vertices $w$ in
layer $i-1$. We call the least one (in the above breadth-first search)
the \emph{parent} of $v$ and $v$ a child of $w$.
Notice that one can draw these layers on
circles. Hence we call such an ordering \emph{spiral like ordering},
see \citep{Pruss:1998a}.

\begin{definition}[BFD-ordering]
  Let $G(V,E)$ be a connected graph with root $v_0$. Then a
  well-ordering $\prec$ of the vertices is called \emph{breadth-first
  search ordering with decreasing degrees} (\emph{BFD}-ordering for
  short) if the following holds for all vertices $v, w\in V$:
  \begin{enumerate}[(B1)]
  \item if $w_1\prec w_2$ then $v_1\prec v_2$ for all children $v_1$
    of $w_1$ and $v_2$ of $w_2$, resp.;
  \item if $v\prec u$, then $d(v)\geq d(u)$.
  \end{enumerate}
  We call a connected graph that has a BFD-ordering of its vertices a
  \emph{BFD-graph}.
\end{definition}

Every graph has for each of its vertices $v$ an ordering with root
$v$ that satisfies (B1). This can be found by a breadth-first
search as described above. However, not all graphs have an ordering
that satisfies (B2); consider the complete bipartite graph $K_{2,3}$.

\begin{thm}
  \label{thm:Cdegseq}
  Let $G$ have greatest maximum eigenvalue in class
  $\mathcal{C}_\pi$. Then there exists a BFD-ordering of $V(G)$ that is 
  consistent with its Perron vector $f$ in such a way that
  $f(u)>f(v)$ implies $u\prec v$ and hence $d(u)\geq d(v)$.
\end{thm}

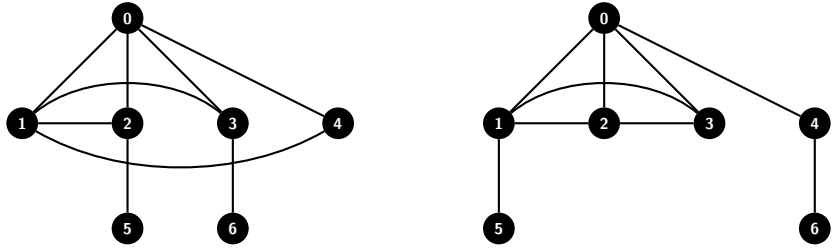
\begin{figure}[ht]
  \centering
  {
\psset{unit=7mm}
\begin{tabular}{c@{\hspace{2cm}}c}
\begin{pspicture}(0,-2.2)(6,2.2)
  \lbvertex{0}{2,2}
  \lbvertex{1}{0,0}
  \lbvertex{2}{2,0}
  \lbvertex{3}{4,0}
  \lbvertex{4}{6,0}
  \lbvertex{5}{2,-2}
  \lbvertex{6}{4,-2}
  \edge{0}{1}
  \edge{0}{2}
  \edge{0}{3}
  \edge{0}{4}
  \edge{1}{2}
  \nccurve[nodesep=0pt,angleA=40,angleB=140]{-}{1}{3}
  \nccurve[nodesep=0pt,angleA=-30,angleB=-150]{-}{1}{4}
  \edge{2}{5}
  \edge{3}{6}
\end{pspicture}
&
\begin{pspicture}(8,-2.2)(14,2.2)
  \lbvertex{0}{10,2}
  \lbvertex{1}{8,0}
  \lbvertex{2}{10,0}
  \lbvertex{3}{12,0}
  \lbvertex{4}{14,0}
  \lbvertex{5}{8,-2}
  \lbvertex{6}{14,-2}
  \edge{0}{1}
  \edge{0}{2}
  \edge{0}{3}
  \edge{0}{4}
  \edge{1}{2}
  \nccurve[nodesep=0pt,angleA=40,angleB=140]{-}{1}{3}
  \edge{1}{5}
  \edge{2}{3}
  \edge{4}{6}
\end{pspicture}
\end{tabular}
}
  \caption{Two BFD-graphs with degree sequence $\pi=(4,4,3,3,2,1,1)$
    that satisfy the conditions of Theorem~\ref{thm:Cdegseq}.\hfill\break
    l.h.s.: $\lambda=3.0918$,
    $f=(0.5291,0.5291,0.3823,0.3823,0.3423,0.1236,0.1236)$,\hfill\break
    r.h.s.:  $\lambda=3.1732$,
    $f=(0.5068,0.5023,0.4643,0.4643,0.1773,0.1583,0.0559)$
    }
  \label{fig:counter-sufficient}
\end{figure}

It is important to note that this condition is not sufficient in
general. Let $\pi=(4,4,3,3,2,1,1)$, then there exist two BFD-graphs
but only one has greatest maximum eigenvalue,
see Figure~\ref{fig:counter-sufficient}. 

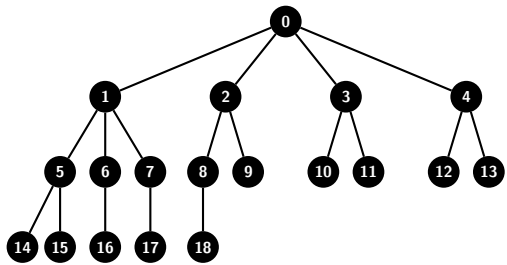
\begin{figure}[ht]
  \centering
  {
  \psset{unit=10mm}
  \begin{pspicture}(-3.5,0)(2.4,-3)
    \lbvertex{0}{0,0}
    \lbvertex{1}{-2.4,-1}
    \lbvertex{2}{-0.8,-1}
    \lbvertex{3}{0.8,-1}
    \lbvertex{4}{2.4,-1}
    \lbvertex{5}{-3,-2}
    \lbvertex{6}{-2.4,-2}
    \lbvertex{7}{-1.8,-2}
    \lbvertex{8}{-1.1,-2}
    \lbvertex{9}{-0.5,-2}
    \lbvertex{10}{0.5,-2}
    \lbvertex{11}{1.1,-2}
    \lbvertex{12}{2.1,-2}
    \lbvertex{13}{2.7,-2}
    \lbvertex{14}{-3.5,-3}
    \lbvertex{15}{-3,-3}
    \lbvertex{16}{-2.4,-3}
    \lbvertex{17}{-1.8,-3}
    \lbvertex{18}{-1.1,-3}
    \edge{0}{1}
    \edge{0}{2}
    \edge{0}{3}
    \edge{0}{4}
    \edge{1}{5}
    \edge{1}{6}
    \edge{1}{7}
    \edge{2}{8}
    \edge{2}{9}
    \edge{3}{10}
    \edge{3}{11}
    \edge{4}{12}
    \edge{4}{13}
    \edge{5}{14}
    \edge{5}{15}
    \edge{6}{16}
    \edge{7}{17}
    \edge{8}{18}
  \end{pspicture}
}
  \caption{A BFD-tree with degree sequence $\pi=(4^2,3^4,2^3,1^{10})$}
  \label{fig:BFD-tree}
\end{figure}

Trees are of special interest. Hence we are looking at the class
$\mathcal{T}_\pi$ of all trees with given sequence $\pi$.
Notice that sequences $\pi=(d_0,\ldots,d_{n-1})$ is a degree 
sequence of a tree if and only if every $d_i>0$ and 
$\sum_{i=0}^{n-1} d_i=2\,(n-1)$, see \citep{Edmonds:1964a}.
In this class there is a single graph with BFD-ordering,
see Figure~\ref{fig:BFD-tree}.

\begin{thm}
  \label{thm:Tdegseq}
  A tree $G$ with degree sequence $\pi$ has greatest maximum
  eigenvalue in class $\mathcal{T}_\pi$ if and only if it is a
  BFD-tree. $G$ is then uniquely determined up to isomorphism.
  The BFD-ordering is consistent with the Perron vector $f$ of $G$ in
  such a way that $f(u)>f(v)$ implies $u\prec v$.
\end{thm}

For a tree with degree sequence $\pi$ a sharp upper bound on the
largest eigenvalue can be found by computing the corresponding 
BFD-tree. Obviously finding this tree can be done 
in $O(n)$ time if the degree sequence is sorted.

We define a partial ordering on degree sequences as follows:
for two sequences $\pi=(d_0,\ldots,d_{n-1})$ and
$\pi'=(d'_0,\ldots,d'_{n-1})$, $\pi\not=\pi'$, we write 
$\pi\lhd\pi'$ if and only if 
$\sum_{i=0}^j d_i\leq \sum_{i=0}^j d'_i$
for all $j=0,\ldots n-1$ (recall that the degree sequences are
non-increasing).
Such an ordering is sometimes called \emph{majorization}.
\begin{thm}
  \label{thm:monotone}
  Let $\pi$ and $\pi'$ two distinct degree sequences of trees with
  $\pi\lhd\pi'$. Let $G$ and $G'$ be trees with greatest maximum
  eigenvalues in classes $\mathcal{C}_\pi$ and $\mathcal{C}_{\pi'}$,
  resp. 
  Then $\lambda(G)<\lambda(G')$.
\end{thm}

We get the following well-known result as an immediate corollary.
\begin{cor}
  \label{thm:Tstar}
  A tree $G$ has greatest maximum eigenvalue in the class of all trees
  with $n$ vertices and $k$ leaves if and only if it is a star with paths 
  of almost the same lengths attached to each of its $k$ leaves.
\end{cor}
\begin{proof}
  The tree sequence $\pi^\ast=(k,2,\ldots,2,1,\ldots,1)$ is maximal
  the class of trees with $k$ pendant vertices w.r.t.\ ordering $\lhd$.
  Thus the statement immediately follows from
  Theorems~\ref{thm:Tdegseq} and \ref{thm:monotone}.
\end{proof}


\section{Proof of the Theorems}
\label{sec:proofs}

We recall that $\lambda(G)$ denotes the maximum eigenvalue of $G$.
Let $N_f(v)=\sum_{uv\in E} f(u)$. Thus the adjacency matrix $A(G)$ can
be defined by $(Af)(v) = N_f(v)$. The Rayleigh quotient of the adjacency
matrix $A(G)$ on vectors $f$ on $V$ is the fraction 
\begin{equation}
  \label{eq:rayleigh}
  \mathcal{R}_G(f)
  = \frac{\langle A f,f\rangle}{\langle f,f \rangle}
  = \frac{\sum_{v\in V} f(v) \sum_{uv\in E} f(u)}{\sum_{v\in V} f(v)^2}
  = \frac{2\sum_{uv\in E} f(u)f(v)}{\sum_{v\in V} f(v)^2}\;.
\end{equation}
By the Rayleigh-Ritz Theorem we find the following well-known
property for the spectral radius of $G$.
\begin{prop}[\citep{Horn:1990a}]
  \label{prop:Courant-Fisher}
  Let $\mathcal{S}$ denote the set of unit vectors on $V$.
  Then 
  \begin{equation*}
    \lambda(G) 
    = \max_{f\in\mathcal{S}} \mathcal{R}_G(f)
    = 2\,\max_{f\in\mathcal{S}} \sum_{uv\in E} f(u)f(v)\;.
  \end{equation*}
  Moreover, if $\mathcal{R}_G(f)=\lambda(G)$ for a (positive)
  function $f\in\mathcal{S}$, then $f$ is an eigenvector
  corresponding to the largest eigenvalue $\lambda(G)$ of $A(G)$,
  i.e., it is a Perron vector.
\end{prop}

The following technical lemma will be useful.
\begin{lem} 
  \label{lem:perron-leaf}
  Let $f$ be the Perron vector of a connected graph $G$. 
  Then $f(u)\geq f(v)$ if and only if $N_f(u)\geq N_f(v)$ .
  Moreover, for each edge $uv\in E$ where $v$ is a pendant vertex and
  $u$ is not, $\lambda(G) = f(u)/f(v)$ and $f(u)>f(v)$.
\end{lem}
\begin{proof}  
  The first statement immediately follows from the positivity of the
  Perron vector and the fact that $f(v) = N_f(v)/\lambda$.
  For the second statement notice that the largest eigenvalue of a
  path with one interior vertex is $\sqrt{2}$.
  Thus the result follows by the well-known fact that 
  $\lambda(H)\leq \lambda(G)$ for a connected subgraph $H$ of $G$.
\end{proof}

The main techniques for proving our theorems is \emph{rearranging} of
edges. We need two standard types of rearrangement steps that we call
\emph{switching} and \emph{shifting}, respectively, in the following.

\begin{lem}[Switching \normalfont{\citep{Rowlinson:1991a,Leydold:1997a}}]
  \label{lem:switching}
  Let $G(V,E)$ be a graph in class $\mathcal{C}_\pi$ with some edges 
  $v_1u_1$ and $v_2u_2$. Assume that $v_1v_2, u_1u_2\notin E$.
  Then we get a new graph $G'(V,E')$ with the same degree sequence
  $\pi$ by replacing $v_1u_1$ and $v_2u_2$ with edges $v_1v_2$ and
  $u_1u_2$ (\emph{switching}).
  Let $f$ is a Perron vector of $G$ then we find
  $\lambda(G')\geq\lambda(G)$, whenever $f(v_1)\geq f(u_2)$ and
  $f(v_2)\geq f(u_1)$.
  The inequality is strict if and only if at least one of these two
  inequalities is strict.
\end{lem}
\begin{proof}
  By removing and inserting edges we obtain
  \[
  \begin{split}
    \mathcal{R}_{G'}(f) - \mathcal{R}_{G}(f) 
    &= 
    \langle A(G') f,f\rangle - \langle A(G) f,f\rangle \\
    &=
    2\left(\sum_{xy\in E'\setminus E}f(x)f(y)
      -  \sum_{uv\in E\setminus E'}f(u)f(v)\right) \\
    &=
    2\,(f(v_1)f(v_2) + f(u_1)f(u_2) - f(v_1)f(u_1) + f(v_2)f(u_2)) \\
    &=
    2\,(f(v_1)-f(u_2))\cdot (f(v_2)-f(u_1)) \\
    &\geq 0\;,
  \end{split}
  \]
  and hence 
  $\lambda(G')\geq\mathcal{R}_{G'}(f)\geq\mathcal{R}_{G}(f)=\lambda(G)$
  by Proposition~\ref{prop:Courant-Fisher}.
  Moreover, $\lambda(G')=\lambda(G)$ if and only if $f$ is also an
  eigenvector corresponding to $\lambda(G')$ on $G'$ and hence
  \[
  \begin{split}
    \lambda(G)f(v_1) &= (A(G)f)(v_1)
    = f(u_1) + \sum_{wv_1\in E \cap E'} f(w) \\
    = \lambda(G')f(v_1) &= (A(G')f)(v_1)
    = f(v_2) + \sum_{wv_1\in E \cap E'} f(w)
  \end{split}
  \]
  and hence $f(u_1) = f(v_2)$. Analogously we find
  $f(v_1)=f(u_2)$.
\end{proof}

\begin{lem}[Shifting 
  \normalfont{\citep{Biyikoglu;Leydold:2006a,Belardo;Marzi;Simic:2006a}}]
  \label{lem:shifting}
  Let $G(V,E)$ be a graph in class $\mathcal{C}_\pi$,
  and let $uv_1\in E$ and $uv_2\notin E$.
  Then we get a new graph $G'(V,E')$ by replacing edge $uv_1$ by the
  edge $uv_2$ (\emph{shifting}).
  Let $f$ is a Perron vector of $G$ then we find
  $\lambda(G')>\lambda(G)$, whenever $f(v_2)\geq f(v_1)$.
\end{lem}
\begin{proof}
  Analogously to the proof of Lemma~\ref{lem:switching} we find
  $\lambda(G')\geq\mathcal{R}_{G'}(f)\geq\mathcal{R}_{G}(f)=\lambda(G)$.
  If equality holded then $f$ would also be a Perron vector of $G'$
  and thus 
  $\lambda(G')f(v_2) 
  = \sum_{xv_2\in E} f(x) + \sum_{yv\in E'\setminus E} f(y) 
  > \sum_{xv_2\in E} f(x) = \lambda(G)f(v_2)$,
  a contradiction.
\end{proof}

\begin{lem}
  \label{lem:shifting-test}
  Let $f$ be the Perron vector of a graph $G$ in $\mathcal{C}_\pi$.
  Let $u$ and $v$ be two vertices with $d(u)>d(v)$. If $f(u)<f(v)$
  then $G$ cannot have greatest maximum eigenvalue in
  $\mathcal{C}_\pi$.
\end{lem}
\begin{proof}
  Let $d(u)-d(v)=c>0$ and assume $f(u)<f(v)$. Then there are (at
  least) $c$ neighbors $w_k$ of $u$ that are not adjacent to $v$.
  When we replace these edges $w_1u,\ldots,w_cu$ by the edges 
  $w_1v,\ldots,w_cv$ we get a new graph $G'$ with the same degree
  sequence $\pi$. 
  The neighbors $c$ can be chosen such that $G'$ remains connected,
  since either $u$ and $v$ have a common neighbor or are adjacent, or
  we can select any of the neighbors of $u$.
  By Lemma~\ref{lem:shifting} we then have $\lambda(G')>\lambda(G)$
  and the statement follows. 
\end{proof}

\begin{lem}
  \label{lem:switching-test}
  Let $f$ be the Perron vector of a graph $G$ in $\mathcal{C}_\pi$.
  Let $vu\in E(G)$ and $vx\notin E(G)$ with
  $f(u)<f(x)\leq f(v)$. If $f(v)\geq f(w)$ for all neigbors $w$ of
  $x$, then $G$ cannot have greatest maximum eigenvalue in
  $\mathcal{C}_\pi$.
\end{lem}
\begin{proof}
  Assume that such vertices exist. Construct a new graph $G'(V,E')$
  with the same degree sequence $\pi$ by replacing edges $vu$ and $xw$
  by edges $vx$ and $uw$. Then by Lemma~\ref{lem:switching},
  $\mathcal{R}_{G'}(f)>\mathcal{R}_G(f)$. It remains to show that we
  can choose vertex $w$ such that $G'$ is connected. Then
  $G'\in\mathcal{C}_\pi$ and hence $G$ cannot have the greatest maximum
  eigenvalue.

  First, notice that there must be a neighbor $p$ of $x$ that is not
  adjacent to $u$, since otherwise $N_f(x)=\sum_{wx\in E} f(w)\leq
  \sum_{yu\in E} f(y) = N_f(u)$ and thus by
  Lemma~\ref{lem:perron-leaf}, $f(x)\leq f(u)$, a contradiction to our
  assumptions. 
  Furthermore, $x$ must have at least two neighbors, since otherwise
  we had by Lemma~\ref{lem:perron-leaf} and assumption $f(x)>f(u)$,
  $f(w)=N_f(x)>N_f(u)\geq f(v)$, a contradiction to $f(w)\leq f(v)$.
  Since $G$ is connected there is a simple path $P_{vx}=(v,\ldots,t,x)$
  from $v$ to $x$.
  Then there are four cases:
  \begin{enumerate}[(1)]
  \item If $vu\notin P_{vx}$ and $ut\notin E(G)$, then we set $w=t$.
  \item Else, if $vu\notin P_{vx}$ and $ut\in E(G)$, then we set $w$
    to one of the neighbors of $x$ that are not adjacent to $u$.
  \item Else, if $vu\in P_{vx}$ and all neighbors not equal $t$ are 
    adjacent to $u$. Then $t$ cannot be adjacent to $u$ and we set
    $w=t$.
  \item Else, $vu\in P_{vx}$ and there exists a neighbor $p$ of $x$,
    $p\not=t$, with $up\notin E(G)$. Then we set $w=p$.
  \end{enumerate}
  In either case $G'$ remains connected. Thus the statement follows.
\end{proof}

\begin{proof}[Proof of Theorem~\ref{thm:Cdegseq}]
  Assume that $G(V,E)$ has greatest maximum eigenvalue in class
  $\mathcal{C}_\pi$. Let $f$ be a Perron vector of $G$.
  Create an ordering $\prec$ by breadth-first search as follows:
  Choose the maximum of $f$ as root $v_0$ in layer $0$;
  append all neighbors $v_1,\ldots,v_{d(v_0)}$ of $v_0$ to the list
  ordered list; these neighbors are ordered such that
  $u\prec v$ whenever $d(u)>d(v)$, or $d(u)=d(v)$ and $f(u)>f(v)$ (in
  the remaining case the ordering can be arbitrary);
  then continue recursively with all vertices $v_1,v_2,\ldots$ 
  until all vertices of $G$ are processed.
  Notice that (B1) holds for this ordering.
  \\
  We first show that $u\prec v$ implies $f(u)\geq f(v)$ for all
  $u,v\in V$.
  Suppose there exist two vertices $v_i$ and $v_j$ with 
  $v_i\prec v_j$ but $f(v_i)<f(v_j)$. 
  Notice that $v_i$ cannot be root $v_0$.
  Let $w_i$ and $w_j$ be the parents of $v_i$ and $v_j$,
  respectively. By construction there are two
  cases: (i) $w_i=w_j$, or (ii) $w_i\prec w_j$. 
  For case (i) we have $d(v_i)\geq d(v_j)$ by construction and 
  $d(v_i)\leq d(v_j)$ by Lemma~\ref{lem:shifting-test} and thus
  $d(v_i)=d(v_j)$. But then we had $v_i\succ v_j$ by the definition of
  our ordering, since $f(v_i)<f(v_j)$, a contradiction.
  \\
  For case (ii) assume that $v_j$ is maximal, i.e., for any other
  vertex $u$ with this property we have $f(u)\leq f(v_j)$.
  Let $v_i$ ($\prec v_j$) be the first vertex (in the ordering of
  $\prec$) with $f(v_i)<f(v_j)$. Hence $f(u)\geq f(v_j)$ for each
  $u\prec v_i$ and we find $f(w_i)\geq f(v_j) > f(v_i)$. 
  Note that $v_j$ cannot be adjacent neither to $w_i$ nor to $v_0$ as
  we then had case (i). Thus $f(w_i)\geq f(u_j)$ for all neighbors
  $u_j$ of $v_j$, since otherwise $v_j$ were not maximal.
  Hence $G$ can not have greatest maximum
  eigenvalue by Lemma~\ref{lem:switching-test}, a contradiction. 
  At last we have to show Property (B2). However, this follows
  immediately from Lemma~\ref{lem:shifting-test}.
\end{proof}

\begin{proof}[Proof of Theorem~\ref{thm:Tdegseq}]
  The necessity condition is an immediate corollary of
  Theorem~\ref{thm:Cdegseq}.
  To show that two BFD-trees $G$ and $G'$ in class
  $\mathcal{T}_\pi$ are isomorphic we use a function $\phi$ that maps
  the vertex $v_i$ in the $i$-th position in the BFD-ordering of
  $G$ to the vertex $w_i$ in the $i$-th position in the
  BFD-ordering of $G'$. By the properties (B1) and (B2) $\phi$ is an
  isomorphism, as $v_i$ and $w_i$ have the same degree and the images
  of neighbors of $v_i$ in the next layer are exactly the neigbors of
  $w_i$ in the next layer. The latter can be seen by looking on all
  vertices of $G$ in the reverse BFD-ordering.
  Thus the proposition follows.
\end{proof}

\begin{proof}[Proof of Theorem~\ref{thm:monotone}]
  Let $\pi=(d_0,\ldots,d_{n-1})$ and $\pi'=(d'_0,\ldots,d'_{n-1})$ be
  two non-increasing tree sequences with $\pi\lhd\pi'$, i.e., 
  $\pi\not=\pi'$, $\sum_{i=0}^j d_i\leq \sum_{i=0}^j d'_i$, and
  $\sum_{i=0}^{n-1} d_i = \sum_{i=0}^{n-1} d'_i = 2(n-1)$.
  Let $G$ have greatest maximum eigenvalue in $\mathcal{T}_\pi$.
  By Theorem~\ref{thm:Tdegseq} $G$ has a BFD-ordering that is
  consistent with $f$, i.e., $f(u)>f(v)$ implies $u\prec v$.
  \\
  First assume that $\pi$ and $\pi'$ differ only in two positions $k$
  and $l$ with $d'_k=d_k+1$ and $d'_l=d_l-1$ (and hence $k<l$ and 
  $d_k\geq d_l>1$).
  Let $v_k$ and $v_l$ be the corresponding vertices in $G$.
  Without loss of generality we assume that $f(v_k)\geq
  f(v_l)$. 
  Since $G$ is a tree and $d(v_l)\geq 2$, there exists a neighbor $w$
  of $v_l$ in layer $h(v_l)+1$ that is not adjacent to $v_k$.
  Thus we can shift edge $v_lw$ by $v_kw$ and get a new
  tree $G'$ with degree sequence $\pi'$ and
  $\lambda(G')>\lambda(G)$ by Lemma~\ref{lem:shifting}. 
  \\
  For two tree sequences $\pi\lhd\pi'$ we can find a sequence of tree
  sequences $\pi=\pi_0\lhd\pi_1\lhd \dots \lhd\pi_k=\pi'$ 
  where $\pi_{i-1}$ and $\pi_i$ ($i=1,\ldots,n$) differ only in two
  positions as described above by the following recursive procedure. 
  For $\pi_{i-1}$ let $j$ be the first position in which $\pi_{i-1}$
  and $\pi'$ differ. Then $d^{(i-1)}_{j}<d'_j$ and we construct
  $\pi_i = (d^{(i)}_0,\ldots,d^{(i)}_{n-1})$ by 
  $d^{(i)}_{j} = d^{(i-1)}_{j}+1$,
  $d^{(i)}_{j+1} = d^{(i-1)}_{j+1}-1$, and  
  $d^{(i)}_{l} = d^{(i-1)}_{l}$ otherwise. If necessary, $\pi_i$ is
  then sorted nonincreasingly. Thus $\pi_i$ again is a tree sequence
  and the statement follows.
\end{proof}


\section{Remarks}
\label{sec:remarks}

In general, we can ask the same questions for Perron vectors of
generalized graph Laplacians, i.e., symmetric matrices with
non-positive off-diagonal entries.
In this paper we showed that switching and shifting operations are
compatible with respect to degree sequences and we used them to find
trees or connected graphs with greatest maximum eigenvalue of the
adjacency matrix. In \citep{Biyikoglu;Leydold:2006a} these operations
were applied to construct graphs with the smallest first eigenvalue of
the so called Dirichlet matrix. Here the corresponding minization
problems are called \emph{Faber-Krahn}-type inequalities. We refer the
interested reader to \citep{Biyikoglu;Leydold;Stadler:2006a} and the
references given therein.

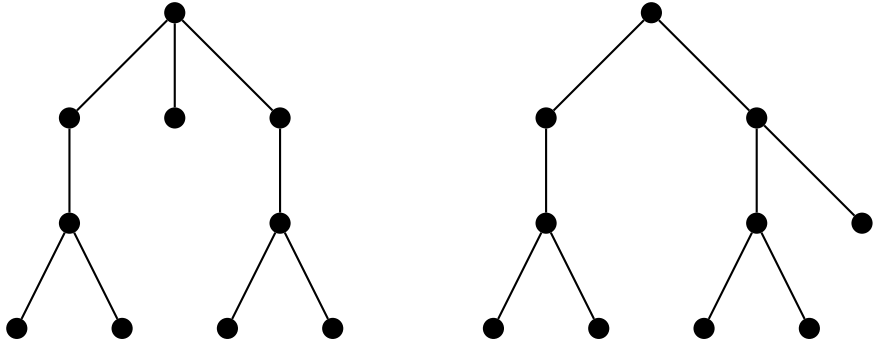
\begin{figure}[ht]
  \centering
  {
\psset{unit=7mm}
\begin{tabular}{c@{\hspace{2cm}}c}
\begin{pspicture}(0,-0.2)(6,6.2)
  \vertex{1}{3,6}
  \vertex{2}{1,4}
  \vertex{3}{3,4}
  \vertex{4}{5,4}
  \vertex{5}{1,2}
  \vertex{6}{5,2}
  \vertex{7}{0,0}
  \vertex{8}{2,0}
  \vertex{9}{4,0}
  \vertex{10}{6,0}
  \edge{1}{2}
  \edge{1}{3}
  \edge{1}{4}
  \edge{2}{5}
  \edge{4}{6}
  \edge{5}{7}
  \edge{5}{8}
  \edge{6}{9}
  \edge{6}{10}
\end{pspicture}
&
\begin{pspicture}(8,-0.2)(15,6.2)
  \vertex{a1}{11,6}
  \vertex{a2}{9,4}
  \vertex{a3}{13,4}
  \vertex{a4}{9,2}
  \vertex{a5}{13,2}
  \vertex{a6}{15,2}
  \vertex{a7}{8,0}
  \vertex{a8}{10,0}
  \vertex{a9}{12,0}
  \vertex{a10}{14,0}
  \edge{a1}{a2}
  \edge{a1}{a3}
  \edge{a2}{a4}
  \edge{a3}{a5}
  \edge{a3}{a6}
  \edge{a4}{a7}
  \edge{a4}{a8}
  \edge{a5}{a9}
  \edge{a5}{a10}
\end{pspicture}
\end{tabular}
}
  \caption{Two trees with degree sequence $(2,2,3,3,3,1,1,1,1,1)$.
    The tree on the l.h.s.\ has smallest maximum eigenvalue
    ($\lambda=2.1010$) among all trees in $\mathcal{C}_\pi$.
    The tree on the r.h.s.\ has a breadth-first ordering of the
    vertices with increasing degree sequences (and thus has lowest
    first Dirichlet eigenvalue). However it does not
    minimize the maximum eigenvalue ($\lambda=2.1067$)
    }
  \label{fig:counter-slo}
\end{figure}

One also might ask whether one can find the smallest maximum
eigenvalue in a class $\mathcal{C}_\pi$ by the same procedure.
It is possible to apply shifting in the proof of
Theorem~\ref{thm:Cdegseq} just the ``other way round''. We then would 
arrive at trees that are constructed by breadth-first search but with
increasing vertex degrees for non-pendant vertices.  
However, this idea does not work.
Figure~\ref{fig:counter-slo} shows a counterexample.

\section*{Acknowledgment}

The authors would like to thank Christian Bey for calling our attention to
eigenvalues of the adjacency matrix of a graph. We thank Gordon
Royle and Brendan McKay for their databases of combinatorial data on
graphs. This was of great help to find the two counterexamples in
Figures~\ref{fig:counter-sufficient} and \ref{fig:counter-slo}.
We also thank the Institute for Bioinformatics of the University in
Leipzig for the hospitality and for providing a scientific working
environment while we wrote down this paper.
The first author is partially supported by
the Belgian Programme on Interuniversity Attraction
Poles, initiated by the Belgian Federal Science Policy Office, and
a grant Action de Recherche Concert\'ee (ARC) of the Communaut\'e
Fran\c{c}aise de Belgique. 



\end{document}